\documentclass[]{amsart}
\usepackage{amsmath,amsfonts,amssymb}
\usepackage{verbatim}
\usepackage{enumerate}
\usepackage{url}
\usepackage[latin1]{inputenc}
\usepackage[english]{babel}
\usepackage{tikz}
\def\N{\mathbb N}

\def\Z{\mathbb Z}

\def\K{\mathbb K}
\def\F{\mathbb F}

\newcommand{\T}{\theta}

\theoremstyle{plain}
\newtheorem{theorem}{Theorem}[section]
\newtheorem{lemma}[theorem]{Lemma}

\newtheorem{proposition}[theorem]{Proposition}
\newtheorem{remark}[theorem]{Remark}

\def\proof{{\it Proof: }}
\def\qed{\hfill\hbox{$\square$}}

\theoremstyle{definition}

\author[F.E. Brochero Mart\'{\i}nez]{F. E. Brochero Mart\'{\i}nez}
\author[Lucas Reis]{Lucas Reis}
\address{
Departamento de Matem\'{a}tica\\
Universidade Federal de Minas Gerais\\
UFMG\\
Belo Horizonte, MG\\
 30123-970\\
 Brazil\\
 }
 \email{fbrocher@mat.ufmg.br }\email{lucasreismat@gmail.com}

\title[Elements of  high order in Artin-Schreier extensions]{Elements of  high order in Artin-Schreier extensions of finite
fields $\F_{q}$}

\date{\today
}

\subjclass[2000]{ }


\begin{document}

\begin{abstract}
In this article,  we find   a lower bound for the order  of the coset  $x+b$  in the Artin-Schreier extension $\F_q[x]/(x^p-x-a) $, 
where $b\in \F_q$ satisfies a  generic  special condition.
\end{abstract}

\maketitle

\section{Introduction}
For many important applications (for example, see \cite{AKS}), it is interesting to find an element of very high order in a finite extension field $\F_{q^n}$. Ideally, one would choose a primitive element, but actually finding such an element is a notoriously hard computation problem.  In fact,  in order to verify that an element is primitive, we need to know the factorization of the integer $q^n-1$ or to solve the discrete logarithms problem in $\F_{q^ n}$.  Now, with the tools currently known these two problem are very hard and   they are the basis of modern cryptography.

On one hand, there are several methods  used to find a small set of elements of $\F_{q^n}$ with at least one primitive  element:  In \cite{Sho}, assuming the extended Riemann hypothesis (ERH),  Shoup has showed  a deterministic polynomial-time search procedure in order to find a primitive element of $\F_{p^2}$; 
Also using ERH, Bach \cite{Bac} gives an  efficiently algorithm in order to  construct a set of $O((\log p)^4/(\log\log p )^3)$ elements,  that contain at least one generator of $\F_p^*$; 
In \cite{Gao}, Gao has given an algorithm to construct high order elements for almost all  
extensions $\F_{q^n}$ of finite fields  $\F_q$,  being   the lower bound no less than $n^{\frac{\log_q n }{ 4 \log_ q ( 2 \log_q n )}- \frac1  2}$. 
Chen \cite{Che} showed how to find, in polynomial time in $N$, an integer $n$ in the interval $[N,2qN]$ and an element $\alpha\in \F_{q^ n}$ with order greater than $5.8^ {n/\log_q n}$.

On the order hand, many works have been done in order to find elements for which a reasonably large lower bound of the order can be guaranteed:  Ahmadi, Shparlinski and  Voloch  \cite{ASV} showed that if $\theta\in\F_{q^{2n}}$ is a primitive $r$-th root of
the unity, where $r=2n+1$ is a prime,  then the Gauss period 
$
\alpha=\theta+\theta^{-1} 
$  has order exceeding  $\exp\left(\left(\pi\sqrt{\frac{2\left(p-1\right)}{3p}}+o\left(1\right)\right)\sqrt{n}\right)$,
 where $p$ is a characteristic of the field (for other works about the order of Gauss period, see \cite{GaSh1,GaSh2}).
  Popovych \cite{Pop1, Pop3} improved
the previous bound and gave a  lower bound for elements of
the more general forms $\theta^e(\theta^f + a)$, ${(\theta^{-f}+a)(\theta^f + a})$ and
$\theta^{-2e}(\theta^{-f}+a)(\theta^f + a)^{-1}$, where $a\in \F_q^*$. In particular, he  proved that the  multiplicative order of the Gauss period $\beta=\theta+\theta^ {-1}$ is not less than $5^ {\sqrt{(r-2)/2}-2}$, for all  $p\ge 5$. 

Finally, Popovych \cite{Pop4} considered the Artin-Schreier extension $\F_{p^p}$ of finite field $\F_p$ and found that an element of order larger than $4^p$ by an elementary method. We emphasize that this  Popovych's result is weaker than  the one,  point out for Shparlinski  in  Voloch's article \cite{Vol}, where, he say  that  the order of any root $t$ of $x^p-x-1$ in $\F_{p^p}$ exceeds $2^{2.54p}\approx 5.81589^p$. Unfortunately, that article does not contain the proof of that limitation and making the computational calculations of that bound, using \emph{Sage Mathematics Software}, we verify that it is true only in the case that $p> 4647$.

In this article, we consider the situation where $x^p-x-a$ is an irreducible polynomial of $\F_q[x]$,  where $p$ is a characteristic of  $\F_q$ and  $a\in \F_p$. We find a lower  bound for the multiplicative order of an  element  of the form $(\theta+ b)$,
 where $\theta$  represents the coset of $x$ in the Artin-Schreier extension $\F_q[x]/(x^p-x-a)$ and $b$ satisfies an special condition.  
We also prove that, the probability that an element of $\F_q$  satisfies such special condition  is  close to $1$ when $q$ is large enough.

Finally,  in the case $q=p$,  we show a lower bound which improves the result obtained by Popovych, but our lower bound does not reach  the one appointed by  Shparlinski-Voloch.

\section{Preliminaries}

Throughout this paper, $\F_q$ denotes a finite field of order $q$, where $q=p^n$ is a power
of an odd prime $p$.   

For each irreducible polynomial $f(x)\in \F_q[x]$, it is known that $\F_q[x]/(f)$ is a finite field with $q^{d}$ elements, where $d=\deg(f)$. Reciprocally, every vector field $\F_{q^d}$ is isomorphic to $\F_q[x]/(f)$ with $f$ an irreducible polynomial of degree $d$.  

There are few known results to ensure the irreducibility of polynomials in a finite field.  For example, Theorems 2.47  and 3.75 in \cite{LiNi} show the necessary and sufficient  conditions for  irreducibility of  cyclotomic polynomials $\Phi_r(x)$ and  binomials  $x^ t-a$, respectively. 
Other well-known result  about the irreducibility of  other family of polynomials is the following. 
\begin{lemma}\label{irredutivel}
The polynomial $x^p-x-a\in \F_q[x]$ is irreducible, if and only if,  it has no roots in $\F_q$.
\end{lemma}  
For the proof of this result, see (Theorem 3.78, \cite{LiNi}).
In particular we have that
\begin{proposition}
Let $n$ be positive integer  and $a\in \F_p^{*}$. The polynomial $f(x)=x^p-x-a$ is irreducible in $\F_q[x]$, if and only if, such that $p\nmid n$.
\end{proposition}

\proof  By Theorem 2.25
in \cite{LiNi}, it is known  that $a=b^p-b$, for some $b\in \F_q$, if and only if, $Tr_{\F_q|\F_p}(a)=a+a^p+\cdots+a^{p^{n-1}}\ne 0$. Since $a\in \F_p$, it follows that $Tr_{\F_q|\F_p}(a)=na$. But   $Tr_{\F_q|\F_p}(a)\ne 0$, if and only if, $p$ does not divide $n$.\qed
\vspace{5mm}

The main results of this paper is the following one:
\begin{theorem}\label{principal} Let $x^p-x-a$ be an irreducible polynomial of $\F_{q}$, with $q=p^n$ ($n\ge 2$) and  $a\in \F_p$. If   $\theta$ is the coset of $x$ in the Artin-Schreier extension $\F_q[x]/(x^p-x-a)$ and  $b\in \F_q$ satisfies that $b \notin \F_{p^m}$, for all $m$ proper divisor of $ n$, then the multiplicative order of $\theta+b$ is lower  bounded by  
$$\frac 1{\pi(p-1)} \sqrt{\frac{2n+1}{2n-1}} \left(\frac {(2n+1)^{2n+1}}{(2n-1)^{2n-1}}\right)^{(p-1)/2} \exp\left(-\frac 1{3(p-1)}\Bigl(\frac {4n^2}{4n^2-1}\Bigr)\right).$$
In particular, for every $\epsilon>0$ and $n> N_\epsilon$,
$$|\langle \theta+b\rangle|> \frac{1}{\pi p} ((e-\epsilon)(2n+1))^{p-1}.$$ 
 
\end{theorem}

And for the case $p=q$, i.e., $n=1$, we obtain

\begin{theorem}\label{Popovych} Let $a\ne 0$ and $b$ be arbitrary  elements of $\F_p$. Then the multiplicative order of $(\theta+b)$ in $\dfrac {\F_p[x]}{(x^p-x-a)}$ is lower bounded by $\frac {\sqrt 3}{\pi p} e^{-\frac 1{12}}\left(\frac {16}3\right)^p$.
\end{theorem}

Observe that using the fields isomorphism $$\begin{matrix}\tau:& \frac{\F_q[x]}{(x^p-x-a)}&\to& \frac{\F_q[x]}{(x^p-x-1)}\\ &h(x)&\mapsto& h(ax)\end{matrix}$$ we only need to prove the Theorem  in the case $a=1$.


\section{The finite field  ${\F_q[x]}/{(x^p-x-1)}$}
Throughout this section, $x^p-x-1$  is an irreducible polynomial of $\F_q[x]$, where
$q=p^n$,   gcd$(p,n)=1$.  Also, $\T$ represents  the coset of $x$ in the Artin-Schreier extensions  $\K:=\F_q[x]/(x^p-x-1)$ and $b\in \F_q\setminus \mathcal 
A_{n}$, where
$$\mathcal A_{n} =\bigcup_{m|n\atop m\ne n} \F_{p^m}.$$

Before we estimate the order of $\theta+b$,  let us show that almost all element of $\F_q$ satisfies  the condition  that we are imposing on $b$.
\begin{theorem}
The number of elements of  $\F_q\setminus \mathcal A_{n}$ is
$\displaystyle\sum_{d|n}p^d\mu(n/d),
$
where $\mu$ is the M\"obius function. In particular, the probability that a chosen element in $\F_q$ does not belong to  $\mathcal A_n$ is greater than $1-\frac {\log_r n}{q^{1-1/r}}$, where $r$ is the smallest prime divisor of $n$.  
\end{theorem}

\proof 
For each positive integer $m$, let $g: \N^*\to \N$  be the function defined by
$$g(m)=|\F_{p^m}\setminus \mathcal A_{m}|.$$  
Clearly,  for each positive integer $m$,  $g(m)$ counts the number of elements in $\F_{p^m}$,  which are  not in any proper subfield of   $\F_{p^m}$. Since each  proper field is of the form $\F_{p^l}$, where $l|m$,  then $$\sum_{d|m} g(d)= |\F_{p^m}|=p^m.$$ 
By the M\"obius Inversion Formula, it follows that 
$$g(m)=\sum_{d|m}p^d\mu(m/d).$$
Now,  let us calculate an upper bound for the number of elements in $\mathcal A_n$. 
Let us suppose that  $p_1^{\alpha_1}\dots p_s^{\alpha_s}$ is the factorization of $n$ in prime factors, where $p_1<\cdots<p_s$. For each  proper divisor $d$ of $n$, there exists a prime $p_i$  $(1\le i\le s)$, such that  $d|(n/p_i)$. 
Thus $\mathcal A_n\subset \bigcup\limits_{1\le i\le s} \F_{p^{n_i}}$ where $n_i=\frac n{p_i}$. 
In particular 
$$|\mathcal A_n|\le \Bigl|\bigcup_{1\le i\le s} \F_{p^{n_i}}\Bigr|\le \sum_{1\le i\le s}p^{n_i}\le p^{n/p_1}\log_{p_1}n= q^{1/p_1}\log_{p_1} n.$$
Therefore,  the probability that a chosen element in $\F_q$ does not belong to  $\mathcal A_n$ is greater than
$$1-\frac {|\mathcal A_n|}{q}\ge  1-\frac {\log_{p_1}(n)}{q^{(1-1/p_1)}}.$$
\qed

This theorem proves that almost all element in $\F_q$ satisfies the condition that we imposed on $b$. 
Now, we need the following technical lemmas:

\begin{lemma}\label{Lemadiferentes}
Let $i$ and $j$ be   integers such that $0\le i, j\le np-1$. If $i\ne j$, then $i+b^{p^i}\ne j+b^{p^j}$.
\end{lemma}
\proof
Let $i_0$ (respectively $j_0$) be the remainder of $i$ (respectively $j$) divided by $n$. We can suppose, without loss of generality, that $i_0\ge j_0$. Clearly, $b^{p^i}=b^{p^{i_0}}$ and $b^{p^j}=b^{p^{j_0}}$. 
Now suppose, by contradiction, that  $i+b^{p^i}=j+b^{p^j}$ and therefore 
\begin{equation}
(j-i)= b^{p^{i_0}}-b^{p^{j_0}}.\label{i_0j_0}
\end{equation}
In the case when  $i_0=j_0$, we have that $j\equiv i\pmod n$, i.e., $j=i+nk$ for some integer $k$ and 
\begin{equation}\label{i_0j_02} 0= b^{p^j}-b^{p^i}=i-j=nk.\end{equation}
It follows that $p$ divides $k$, what is impossible because $0<|i-j|<np$.  

Thus  $0<i_0-j_0<n$, and taking the $p^{n-j_0}-$th power in (\ref{i_0j_02}),  we have
$$j-i=b^{p^{n+i_0-j_0}}-b=b^{p^{i_0-j_0}}-b.$$
Thereby, there exists $0\le t< n$ such that $b^{p^t}-b\in \F_p$, or equivalently $b^{p^{t+1}}-b^p=(b^{p^t}-b)^p=b^{p^t}-b$. This last equation can be rewritten as
$b^p-b=(b^p-b)^{p^t}$, i.e., $b^p-b$ is an element of $\F_{p^t}$.
Furthermore, if $b\notin \F_{p^t}$, by Lemma \ref{irredutivel}, the polynomial $x^p-x-(b^p-b)$ is an irreducible polynomial of $\F_{p^t}$. We obtain, in any case, that  $b\in \F_{p^{pt}}$. Since $b$ is also in $\F_{p^n}$, we conclude that $b$ belongs to
$$\F_{p^{pt}}\cap \F_{p^n}=\F_{p^{\text{gcd}(pt,n)}}=\F_{p^{\text{gcd}(t,n)}},$$
where gcd$(t,n)<n$ is a proper divisor of $n$ and so we have a contradiction with the choice of $b\not\in \mathcal A_n$.\qed

\begin{lemma}\label{injetivo}Let $t, s$ be nonnegative integers such that $0\le t+s\le p-1$ and let $I_{s,t}$ be the subset of  $\Z^{np}$  such that $\vec r:=(r_0,r_1,\dots, r_{np-1})\in I_{s,t}$ if and only if
$$\sum_{ 0\le j\le np-1\atop r_j<0}(- r_j)\le t\quad\text{and} \sum_{ 0\le j\le np-1\atop r_j>0} r_j\le s$$
Then the function 
$$\begin{matrix}
\Lambda:&I_{s,t}&\to& G\\
&\vec r&\mapsto&\prod\limits_{0\le j\le np-1}(\theta+b)^{r_j p^j},
\end{matrix}
$$
where $G=\langle \theta+b\rangle\le \K^*$, 
is one to one.

\end{lemma}

\proof
Since  $\T$ is the  coset of $x$ in the quotient field $\K=\frac{\F_q[x]}{x^p-x-1}$, then  each element of $\K$ is the coset of a unique $h(\theta)$, where $h$ is a polynomial in $\F_q[x]$ of degree at most
 $p-1$. 
In addition,   
 $\T^p=\T+1$ and, accordingly, for all $j\in \N$,
 $$\T^{p^{j+1}}=(\T^{p})^{p^j}=(\T+1)^{p^j}=\T^{p^j}+1. $$
It follows,  inductively, that
$$
\T^{p^j}=\T+j \quad \text{for all $j\ge 1$,}\label{potencia}
$$
and, thereby, for each $\vec{r}=(r_0,\dots, r_{np-1})\in I_{s,t}$ 
$$\Lambda(\vec r)=\prod_{0\le i\le np-1}(\T+b)^{r_ip^i}=\prod_{0\le i\le np-1}(\T+i+b^{p^i})^{r_i}.$$
Now, suppose that  $\vec s=(s_0,\dots, s_{np-1})$ is another element of $I_{s,t}$ such that $\Lambda(D)=\Lambda(E)$,  i.e.,
$$\prod\limits_{0\le i\le np-1}(\T+i+b^{p^i})^{r_i}= \prod\limits_{0\le j\le np-1}(\T+j+b^{p^j})^{s_j},$$
thus, the polynomial
$$F(x)=\prod\limits_{0\le i\le np-1\atop r_i>0}(x+i+b^{p^i})^{r_i}\prod\limits_{0\le j\le np-1\atop s_j<0}(x+j+b^{p^j})^{-s_j}$$
is congruent to  the polynomial
$$ G(x)=\prod\limits_{0\le j\le np-1\atop s_j>0}(x+j+b^{p^j})^{s_j}\prod\limits_{0\le i\le np-1\atop r_i<0}(x+i+b^{p^i})^{-r_i}$$
modulo ${x^p-x-1}$. 

 Since $\deg(F)\le s+t\le p-1$ and $\deg(G)\le s+t<p-1$, it follows that $F(x)=G(x)$.
Further, by Lemma \ref{Lemadiferentes}, we know that  $x+i+b^{p^i} \ne x+j+b^{p^j}$,  for  all $0\le i< j\le np-1$,
  therefore $\vec r=\vec s$, as we want to prove.
  \qed
  
We emphasize that, in the last step of  Lemma, is essential the condition that we imposed on $b$.

\begin{lemma}\label{I_st}  Let $I_{s,t}$ be as in the Lemma \ref{injetivo}. Then
\begin{equation}
|I_{s,t}|=\sum_{j=0}^t\sum_{i=0}^s \binom{np}{i} \binom{np-i}{j} \binom{s}{i}\binom{t}j.\label{numeroelementos}
\end{equation}
In particular, 
$$|I_{s,t}|> \binom{np+t-s}t\binom{np+s}s
.$$
\end{lemma}

\proof Observe that, for each $j\le t$ and $i\le s$, we can select $j$ coordinates of $\vec r$ to be negative and $i$ coordinates to be positive and  this  choice can be done of   $\binom{np}{i} \binom{np-i}{j}$ ways.   Besides,  the number of positive  solution of $x_1+x_2+\cdots+x_i\le s$ is  $\binom{s}{i}$
and the number of positive  solution of $x_1+x_2+\cdots+x_j\le t$ is  $\binom{t}{j}$. Thus, for each pair  $i, j$, there exist $\binom{np}{i} \binom{np-i}{j} \binom{s}{i}\binom{t}j$ elements of $I_{s,t}$ and then, adding over all $i$ and $j$, we conclude the equality (\ref{numeroelementos}).
In addition
\begin{align*}
 |I_{s,t}|&\ge \sum_{i=0}^s \binom{s}{i}\binom{np}{i}\Bigl(\sum_{j=0}^t  \binom{np-i}{j}\binom{t}j\Bigr)\\
&=\sum_{i=0}^s \binom{s}{i}\binom{np}{i}\binom{np+t-i}t\\
&> \binom{np+t-s}t\sum_{i=0}^s \binom{s}{i}\binom{np}{i}\\
&=\binom{np+t-s}t\binom{np+s}s
\end{align*}
\qed

 Before proceeding to prove  the main Theorems, we need the following technical Lemma, that is essentially  a good application of Stirling approximation.
\begin{lemma}[\cite{Sas} Corollary 1]\label{sasvari} For all $s> 0$ and  $r> 1$, we have
$$c_r\cdot d_r^s\cdot \frac 1{\sqrt s}\cdot \Theta(r,s)< \binom{rs} s <c_r\cdot d_r^s\cdot \frac 1{\sqrt s},$$
where
$$c_r=\sqrt{\frac r{2\pi (r-1)}}, \quad d_r=\frac{r^r}{(r-1)^{r-1}}$$
and
$$\Theta(r,s)=\exp\left(-\frac 1{12s}\left(1+\frac 1{r(r-1)}\right)\right).$$
\end{lemma}

We  emphasize that these upper and lower  bounded  are very close  when $s\gg 0$.

\section{Proof of Theorem \ref{principal}}
By Lemma \ref{injetivo}, we know that 
$|\langle \theta+b\rangle|\ge |I_{s,t}|$,
for all nonnegative  integers $s$ and $t$ such that $s+t\le p-1$.   
So,  by Lemma \ref{I_st}, we have that
\begin{align}
|\langle \theta+b\rangle|&>\max_{0\le s+t\le p-1} \binom{np+t-s}t\binom{np+s}s\nonumber\\
&>\binom{np}{(p-1)/2}\binom{np+(p-1)/2}{(p-1)/2}.\label{binomios}
\end{align}
Now, using Lemma \ref{sasvari}, each binomial coefficient can be bounded by
\begin{align*}
\binom{np}{(p-1)/2}&>\binom{2n(p-1)/2}{(p-1)/2}\\
&>\sqrt{\frac{2n}{\pi(2n-1)(p-1)} } \left(\frac {(2n)^{2n}}{(2n-1)^{2n-1}}\right)^{\frac{p-1}2}  \tilde\Theta(2n-1)
\end{align*}
and
\begin{align*}
\binom{np+(p-1)/2}{(p-1)/2}&>\binom{(2n+1)(p-1)/2}{(p-1)/2}\\
&> \sqrt{\frac{2n+1}{\pi(2n)(p-1)} } \left(\frac {(2n+1)^{2n+1}}{(2n)^{2n}}\right)^{\frac{p-1}2} 
 \tilde\Theta(2n+1)
\end{align*}
where $\tilde\Theta(z)=\exp\left(-\frac 1{6(p-1)}\Bigl(1+\frac {1}{z(z-1)}\Bigr)\right)$.

Multiplying these two inequalities and simplifying, we conclude  that
$$|\langle \theta+b\rangle|>\frac 1{\pi(p-1)} \sqrt{\frac{2n+1}{2n-1}} \left(\frac {(2n+1)^{2n+1}}{(2n-1)^{2n-1}}\right)^{(p-1)/2} \exp\left(-\frac 1{3(p-1)}\Bigl(\frac {4n^2}{4n^2-1}\Bigr)\right).$$
Therefore, we obtain the first part of the Theorem. 

For the second part, observe that  the sequence $\{a_n\}_{n\in \N}$ defined  for each $n\ge 2$, as  $a_n:= \left(\frac {2n+1}{2n-1}\right)^{\frac{(2n-1)(p-1)+1}2}$, is an increasing sequence satisfying  $$a_2>\sqrt{\frac 53}(2.1516)^{p-1}\quad\text{and}\quad \lim\limits_{n\to \infty} a_n=e^{p-1}.$$
 Therefore, for $n\ge 2$, we can find a simpler but weaker estimate
\begin{align*}
|\langle \theta+b\rangle|&>\frac 1{\pi(p-1)} (a_n(2n+1))^{p-1} \exp\left(-\frac 1{3(p-1)}\Bigl(\frac {4n^2}{4n^2-1}\Bigr)\right)\\
&>\frac {\sqrt 5}{\sqrt 3\pi(p-1)}(2.1516(2n+1))^{p-1} \exp\left(-\frac {16}{45(p-1)}\right).
\end{align*}
In the case $n$ large enough, we have that
$$(e-\epsilon)^{p-1}<a_n<e^{p-1} \quad\text{and}\quad \exp\left(-\frac 1{3(p-1)}\Bigl(\frac {4n^2}{4n^2-1}\Bigr)\right)> 1-\frac {16}{45(p-1)},$$ 
therefore 
\begin{align*}|\langle \theta+b\rangle|&>  \frac {45p-61}{45\pi (p-1)^2} ((e-\epsilon)(2n+1))^{p-1}\\
&>\frac {1}{\pi p} ((e-\epsilon)(2n+1))^{p-1},
 \end{align*}
  as we want to prove.
\qed

The following table the lower bounded of $|\langle \theta+b\rangle|$, for some values of $n$, where  the value of $p$ appears as a parameter
\begin{center}
\begin{tabular}{c|c}
$n$& $\pi p \cdot |\langle \theta+b\rangle|$\\ \hline
$2$&$12.22377^ p$\\  \hline
 $ 3$&$ 17.65835 ^ p$\\ \hline 
 $4 $&$  23.09586^ p$\\ \hline 
 $5$&$ 28.53356^ p$\\ \hline 
 $10 $&$ 55.71983 ^ p$\\ \hline 
 $ 100$&$ 545.01494 ^ p$\\ \hline 
 $1000 $&$ 5437.92274 ^ p$\\ \hline 
 $10000 $&$ 54366.9957 ^ p$\\ \hline 
\end{tabular}
\end{center}

\section{Proof of Theorem \ref{Popovych}}
 The polynomial $x^p-x-1$,  is always an irreducible polynomial of $\F_p[x]$ and   the condition imposed on $b$ is empty.  So, 
by Lemma \ref{I_st}, we have that
\begin{align*}
|\langle \theta+b\rangle|&>\max_{s+t= p-1} \binom{p+t-s}t\binom{p+s}s\\
&=\max_{0\le s\le p-1}  \binom{2p-1-2s}{p-1-s}\binom{p+s}s.\\
&=\max_{0\le \lambda\le \frac{p-1}p \atop p \lambda \in \N}  \binom{2p-1-2p\lambda}{p-1-p\lambda}\binom{p+p\lambda}{p\lambda}\\
&=\frac 12\max_{0\le \lambda\le \frac{p-1}p \atop p \lambda \in \N}  \binom{p(2-2\lambda)}{p(1-\lambda)}\binom{p(1+\lambda)}{p\lambda }.\label{binomios2}
\end{align*}
The same way, using  Lemma \ref{sasvari},  we obtain that 
$$|\langle \theta+b\rangle|>\max_{0\le \lambda\le \frac{p-1}p \atop p\lambda \in \N} \frac 1{\pi p} \sqrt{\frac{1+\lambda}{2\lambda(1-\lambda)}} \left(\frac{4^{1-\lambda} (1+\lambda)^{1+\lambda}}{\lambda^\lambda}\right)^p \Theta(2,p(1-\lambda))\Theta\left(\frac {1+\lambda}\lambda, p\lambda\right),
$$
in particular, taking    $\lambda=\frac 13$, it follows that
$$|\langle \theta+b\rangle|> \frac {\sqrt 3}{\pi p} e^{-\frac 1{12}}\left(\frac {16}3\right)^p.\eqno\qed
$$

\begin{remark}
In summary,  for the case $\F_{p^p}$ and $p\gg 0$, we observe that lower bound of $|\langle \theta+b\rangle |$  is $O(4^p)$ in Popovych's paper, 
$O( 5.81589^p)$ in Voloch's article and  $O(5.3333^p)$ in our result. 
\end{remark}

\end{document}